\documentclass[12pt]{amsart}
\usepackage{amsmath,latexsym,amscd,amsbsy,amssymb,amsfonts,amsthm,fleqn,leqno,
euscript, graphicx, texdraw}
\numberwithin{equation}{section}

\newcommand{\cl}{\mathrm{cl}}

\newcommand{\I}{\mathbb I}

\newcommand{\diam}{\mathrm{diam}}

\newtheorem{thm}{Theorem}[section]
\newtheorem{pro}[thm]{Proposition}
\newtheorem{lem}[thm]{Lemma}

\newtheorem*{MainTheorem1}{Theorem \ref{main1}}
\newtheorem*{MainTheorem2}{Theorem \ref{thm2}}
\newtheorem*{MainTheorem3}{Theorem \ref{C}}
\newtheorem*{MainTheorem4}{Theorem~\ref{wid}}
\newtheorem*{MainTheorem5}{Theorem~\ref{inf}}

\newtheorem{qu}[thm]{Question}

\theoremstyle{definition}
\newtheorem{dfn}[thm]{Definition}
\newtheorem{re}[thm]{Remark}

\theoremstyle{remark}
\newtheorem{claim}[thm]{Claim}

\newcommand{\inte}{\mathrm{int}}

\newcommand{\bd}{\mathrm{bd}}

\begin{document}

\title[Generalized Cantor manifolds and homogeneity]
{Generalized Cantor manifolds and homogeneity}

\author{A. Karassev}
\address{Department of Computer Science and Mathematics,
Nipissing University, 100 College Drive, P.O. Box 5002, North Bay,
ON, P1B 8L7, Canada} \email{alexandk@nipissingu.ca}

\author{P. Krupski}
\address{Mathematical Institute, University of Wroc\l aw, pl.
Grunwaldzki 2/4, 50--384 Wroc\l aw, Poland}
\email{krupski@math.uni.wroc.pl}

\author{V. Todorov}
\address{Department of Mathematics, UACG, $1$ H. Smirnenski blvd.,
$1046$ Sofia, Bulgaria} \email{vtt-fte@uacg.bg}

\author{V. Valov}
\address{Department of Computer Science and Mathematics,
Nipissing University, 100 College Drive, P.O. Box 5002, North Bay,
ON, P1B 8L7, Canada} \email{veskov@nipissingu.ca}

\date{\today}
\thanks{The first author was partially supported by NSERC
Grant 257231-04.}

\thanks{The paper originated during the second author's stay at the Nipissing University in 2006--2007}

\thanks{The last author was partially supported by NSERC
Grant 261914-03.}


\begin{abstract}
A classical theorem of Alexandroff states that every $n$-dimensional
compactum $X$ contains an $n$-dimensional Cantor manifold. This
theorem has a number of generalizations obtained by various authors.
We consider extension-dimensional and infinite dimensional analogs
of strong Cantor manifolds, Mazurkiewicz manifolds, and
$V^n$-continua, and prove corresponding versions of the above
theorem. We apply our results to show that each homogeneous
metrizable continuum  which is not in a given  class $\mathcal C$ is
a strong Cantor manifold (or at least a Cantor manifold) with
respect to $\mathcal C$. Here, the class $\mathcal C$ is one of four
classes that are defined in terms of dimension-like invariants.
  A class of
spaces having bases of neighborhoods satisfying certain special
conditions is also considered.
\end{abstract}
\maketitle\markboth{}{}



\tableofcontents


\section{Introduction}
All spaces in this paper are assumed to be at least normal.

Cantor manifolds were introduced by Urysohn \cite{u} as a
generalization of Euclidean manifolds. Recall that a space $X$ is a
{\em Cantor $n$-manifold} if $X$ cannot be separated by a closed
$(n-2)$-dimensional subset. In other words, $X$ cannot be the union
of two proper closed sets whose intersection is of covering
dimension $\leq n-2$. Alexandroff \cite{ps} introduced the stronger
notion of $V^n$-continua: a compactum $X$ is a {\em $V^n$-continuum}
if for every two closed disjoint subsets $X_0$, $X_1$ of $X$, both
having non-empty interior in $X$, there exists an open cover
$\omega$ of $X$ such that there is no partition $P$ in $X$ between
$X_0$ and $X_1$ admitting an $\omega$-map into a space $Y$ with
$\dim Y\leq n-2$. Another specification of Cantor manifolds was
considered by Had\v{z}iivanov \cite{h}: $X$ is a {\em strong Cantor
$n$-manifold} if for arbitrary representation
$X=\bigcup_{i=1}^{\infty}F_i$, where all $F_i$ are proper closed
subsets of $X$, we have $\dim(F_i\cap F_j)\geq n-1$ for some $i\neq
j$.

Obviously, strong Cantor $n$-manifolds are Cantor $n$-manifolds.
Moreover,  every $V^n$-continuum is a strong Cantor $n$-manifold
\cite{ht} and none of the above inclusions is reversible (see
\cite{hh}, \cite{lel} and the Appendix).

In the present paper we generalize these notions by considering a
general dimension function $D_{\mathcal{K}}$ which captures the
covering dimension, cohomological dimension $\dim_G$ with respect to
any Abelian group $G$, as well as the extraordinary dimension
$\dim_L$ with respect to a given $CW$-complex $L$.

More precisely, a sequence $\mathcal{K}=\{K_0,K_1,..\}$ of
$CW$-complexes is called a {\em stratum} for a dimension theory
\cite{dr} if
\begin{itemize}\item
for each space $X$ admitting a perfect map onto a metrizable space,
$K_n\in AE(X)$ implies both $K_{n+1}\in AE(X\times\I)$ and
$K_{n+j}\in AE(X)$ for all $j\geq 0$.
\end{itemize}
Here, $K_n\in AE(X)$ means that $K_n$ is an absolute extensor for
$X$. Given a stratum $\mathcal{K}$, we can define a dimension
function $D_{\mathcal{K}}$ in a standard way:
\begin{enumerate}
\item
$D_{\mathcal{K}}(X)=-1$ iff $X=\emptyset$;
\item $D_{\mathcal{K}}(X)\le n$ if
$K_n\in AE(X)$ for $n\ge 0$; if $D_{\mathcal{K}}(X)\le n$ and
$K_m\not\in AE(X)$ for all $m<n$, then $D_{\mathcal{K}}(X)= n$;
\item
$D_{\mathcal{K}}(X)=\infty$ if $D_{\mathcal{K}}(X)\le n$ is not
satisfied for any $n$.
\end{enumerate}

Since every $CW$-complex $K$ with the weak topology is homotopically
equivalent to $K$ equipped with the metric topology, we can assume
that all $K_i\in\mathcal{K}$ are considered with the metric
topology.

 If $\mathcal{K}=\{\mathbb{S}^0,\mathbb{S}^1,..\}$, we
obtain the covering dimension $\dim$. The stratum
$\mathcal{K}=\{\mathbb{S}^0,K(G,1),..,K(G,n),..\}$,  $K(G,n)$,
$n\geq 1$, being the Eilenberg-MacLane complexes for a given group
$G$, determines the cohomological dimension $\dim_G$. Moreover, if
$L$ is a fixed $CW$-complex and
$\mathcal{K}=\{L,\Sigma(L),..,\Sigma^n(L),..\}$, where $\Sigma^n(L)$
denotes the $n$-th iterated suspension of $L$, we obtain the
extraordinary dimension $\dim_L$ introduced recently by Shchepin
\cite{es:98} and considered in details by Chigogidze \cite{ch:03}.

According to the countable sum theorem in extension theory, it
follows directly from the above definition that
$D_{\mathcal{K}}(X)\leq n$ implies $D_{\mathcal{K}}(A)\leq n$ for
any $F_{\sigma}$-subset $A\subset X$.

Now, it is clear how to define Cantor $n$-manifolds, strong Cantor
$n$-manifolds and $V^n$-continua with respect to $D_{\mathcal{K}}$,
where $\mathcal{K}$ is a fixed stratum. Furthermore, we  consider
quite general concepts of  Mazurkiewicz manifolds, strong Cantor
manifolds and Cantor manifolds with respect to some classes of
finite or infinite-dimensional  spaces. We define them following the
idea and some terminology from \cite{hh, ht}.

A non-empty class of spaces $\mathcal{C}$ is said to be {\em
admissible} if it satisfies the following conditions:
\begin{itemize}
\item[(i)] $\mathcal{C}$ contains all topological copies of any element $X\in\mathcal{C}$ ;
\item[(ii)] if $X\in\mathcal{C}$, then each $F_\sigma$-subset of $X$ belongs to $\mathcal{C}$.
\end{itemize}

\begin{dfn}\label{dfn1}
A space $X$  is a {\em Mazurkiewicz manifold with respect to an
admissible class $\mathcal{C}$} if for every two closed, disjoint
subsets $X_0,X_1\subset X$, both having non-empty interiors in $X$,
and every $F_\sigma$-subset $F\subset X$ with $F\in\mathcal{C}$,
there exists a continuum in $X\setminus F$ joining $X_0$ and $X_1$.
\end{dfn}
The notion of a Mazurkiewicz manifold has its roots in  the
classical Mazurkiewicz theorem saying that no region in the
Euclidean $n$-space can be cut by a  subset of dimension $\leq
n-2$~\cite{re:95}.  Recall that a set $P$ (not necessarily closed)
{\em cuts} a space $X$ between two subsets $X_0$ and $X_1$ of $X$ if
$X_0$, $X_1$, and $P$ are disjoint, and for any continuum $C$ such
that $C\cap X_i\ne\emptyset$, $i=0,1$, we have $C\cap P
\ne\emptyset$; $P$ {\em cuts} $X$ if it cuts $X$ between a pair of
distinct points.

One can easily prove, using Lemma~\ref{mazur1}, that if no
$F_\sigma$-subset from an admissible class $\mathcal{C}$ cuts a
compact space $X$, then $X$ is a Mazurkiewicz manifold with respect
to $\mathcal{C}$; the converse implication holds for locally
connected compact spaces $X$.

\begin{dfn}
A space $X$  is a {\em strong Cantor manifold with respect to an
admissible  class $\mathcal{C}$} if $X$ can not be represented as
the union
\begin{align}\displaystyle \label{eq1}& X=\bigcup_{i=0}^\infty F_i \quad \text{with}\quad \bigcup_{i\ne j}(F_i\cap F_j)\in \mathcal C\\
&\text{where all $F_i$ are proper closed subsets of  $X$}.
\notag\end{align}
 \end{dfn}

\begin{dfn}
A space $X$  is a {\em  Cantor manifold with respect to an
admissible  class $\mathcal{C}$} if $X$ cannot be separated by a
closed subset which belongs to $\mathcal{C}$.
\end{dfn}

Four specifications of $\mathcal{C}$ will be  considered:
\begin{enumerate}
\item
 the class $\mathcal D_{\mathcal K}^k$ of at most
$k$-dimensional spaces with respect to dimension $D_{\mathcal{K}}$,

\item the class $\mathcal D_{\mathcal{K}}^{<\infty}$ of strongly countable
$D_{\mathcal K}$-dimensional spaces, i.e. all spaces represented as
a countable union of closed finite-dimensional subsets with respect
to $D_{\mathcal{K}}$,

\item
the class $\mathbf C$ of paracompact $C$-spaces,

 and
\item the class $\mathcal {WID}$
of weakly infinite-dimensional spaces.
\end{enumerate}

Recall that $X$ is said to be {\em strongly infinite-dimensional} if
for any sequence $\{(A_n,B_n)\}_{n\geq 1}$ of pairs of disjoint
closed sets in $X$ and any sequence of closed partitions $C_n\subset
X$ separating $A_n$ and $B_n$ the intersection $\bigcap_{n\geq
1}C_n$ is non-empty. Spaces which are not strongly
infinite-dimensional are called {\em weakly infinite-dimensional}.

A space $X$ is said to be a {\em $C$-space} (or has {\em property
$C$}) \cite{re:95} if for every sequence $\{\omega_n\}_{n\geq 1}$ of
open covers of $X$ there exists a sequence $\{\gamma_n\}_{n\geq 1}$
of open disjoint families in $X$ such that each $\gamma_n$ refines
$\omega_n$ and $\bigcup_{n\geq 1}\gamma_n$ is a cover of $X$.

Every finite-dimensional paracompact space as well as every
count\-able\--dimensional metrizable space is a $C$-space, but there
exist  metrizable $C$-spaces which are not countable-dimensional
\cite{rp}. Moreover, compact $C$-spaces form a proper subclass of
weakly infinite-dimensional compact spaces \cite{bo}.

Every compact Mazurkiewicz manifold with respect to any admissible
class $\mathcal{C}$ is a strong Cantor manifold with respect to
$\mathcal{C}$ (see Proposition~\ref{mazur}) and strong Cantor
manifolds with respect to $\mathcal{C}$ are  Cantor manifolds with
respect to $\mathcal{C}$.

The following theorems are amongst the main results of the paper.

\begin{MainTheorem1}
Any compact space $X$ with $D_{\mathcal{K}}(X)=n$ contains a closed
subset $M$ such that $D_{\mathcal{K}}(M)=n$ and $M$ is both a
$V^n$-continuum and a Mazurkiewicz manifold with respect to the
class $\mathcal D_{\mathcal{K}}^{n-2}$.
\end{MainTheorem1}

\begin{MainTheorem5} If a compact space $X$ has dimension $D_{\mathcal{K}}(X)=\infty$,
then either $X$ contains closed subsets of arbitrary large finite
dimension $D_{\mathcal{K}}$ or $X$ contains a compact Mazurkiewicz
manifold with respect to the class $\mathcal
D_{\mathcal{K}}^{<\infty}$.
\end{MainTheorem5}

\begin{MainTheorem3} Any compact space without property $C$
contains a closed set which is a Mazurkiewicz manifold with respect
to the class $\mathbf C$.
\end{MainTheorem3}

\begin{MainTheorem4} Any metrizable strongly infinite-dimensional compact
space contains a closed set which is a Mazurkiewicz manifold with
respect to the class $\mathcal {WID}$.
\end{MainTheorem4}

Based on these theorems, we prove the following result.

\begin{MainTheorem2}
Each metrizable homogeneous continuum $X\notin\mathcal{C}$ is a
strong Cantor manifold with respect to class $\mathcal{C}$ provided
that:
 \begin{enumerate}
\item
  $\mathcal{C}$ is any
of the  following three classes:  $\mathcal{WID}$, $\mathbf C$,
$\mathcal D_{\mathcal K}^{n-2}$
 $($in the latter case  we additionally assume
$D_{\mathcal{K}}(X)=n$$)$;\\
or
 \item
$\mathcal C=\mathcal D_{\mathcal K}^{<\infty}$ and   $X$ does not
contain closed subsets of arbitrary large finite dimension
$D_{\mathcal K}$.
\end{enumerate}
 \end{MainTheorem2}

Theorem~\ref{C} is totally new, while some particular weaker cases
of Theorems \ref{main1}, \ref{inf} and \ref{wid} were proved by
different authors. Let us mention the classical result that every
compact space $X$ with the covering dimension $\dim X=n$ contains an
$n$-dimensional Cantor $n$-manifold (with respect to $\dim$)
established independently by Hurewicz-Menger \cite{hm} and Tumarkin
\cite{tu} for metrizable spaces, and by Alexandroff \cite{ps1} for
any compact spaces. For $V^n$-continua with respect to  $\dim$, this
theorem was obtained by Alexandroff \cite{ps} (metrizable compact
spaces) and Kuz'minov \cite{ku} (arbitrary compact spaces). Both
Alexandroff's and Kuz'minov's proofs are based on cohomological
methods. An elementary proof was given by Hamamd\v{z}iev \cite{ham}.
For strong Cantor $n$-manifolds with respect to $\dim_G$, Theorem
\ref{main1} appeared in \cite{hs}.

A classical counterpart  of Theorem~\ref{inf} saying that each
infinite-dimensional compact space $X$  contains either closed
subsets of arbitrary large finite dimension or  a Cantor
$\infty$-manifold $M$ (i.e., no finite-dimensional subset separates
$M$) was proved by Tumarkin~\cite{tu1}.

The fact that each strongly infinite-dimensional compact metric
space contains a compact strongly infinite-dimensional Cantor
manifold $M$ (i.e., no weakly infinite-dimensional closed subset of
$M$ separates $M$) is due to Skljarenko (see~\cite[p. 550]{ap}).

One of the main technical tools in proving Theorem~\ref{main1} is an
extension theorem, see Proposition~\ref{extension1}. In its turn,
Proposition~\ref{extension1} implies another general extension
theorem (Proposition~\ref{extension2}) whose analogues were
established by Holszty\'{n}ski \cite{ho}, Had\v{z}iivanov \cite{h1}
and Dijkstra \cite{di} for covering dimension.
Had\v{z}iivanov-Shchepin \cite[Theorem 1]{hs} also formulated
similar to Proposition~\ref{extension2} statement for cohomological
dimension. However, we were not able to verify some details in their
proof. Instead of following the arguments of the above authors, we
base our proofs of Proposition~\ref{extension1} and
Proposition~\ref{extension2} on a completely different idea.

It was proved in~\cite{kru} that every homogeneous metrizable,
locally compact, connected space $X$ with the covering dimension
$\dim X=n\le\infty$ is a Cantor $n$-manifold; in case where $X$ is
strongly  infinite-dimensional, it is a strongly
infinite-dimensional Cantor manifold. Theorem~\ref{thm2}
significantly generalizes  those  results.

The final section contains examples distinguishing the following
four classes (with respect to $\dim$): Cantor $n$-manifolds, strong
Cantor $n$-manifolds, Mazurkiewicz $n$-manifolds and $V^n$-continua.

\section{Mazurkiewicz $n$-manifolds and $V^n$-continua with respect to dimension $D_{\mathcal{K}}$}

\begin{pro}\label{mazur} Let
$\mathcal{C}$ be an admissible class of spaces. Then  every compact
Mazurkiewicz manifold with respect to $\mathcal{C}$ is a strong
Cantor manifold with respect to $\mathcal{C}$.
\end{pro}

\begin{proof}
Suppose $X$ is a compact Mazurkiewicz manifold with respect to
$\mathcal{C}$ but not a strong Cantor manifold with respect to
$\mathcal{C}$. Then $X=\bigcup_{i\geq 0}F_i$ with  $F_i$ being
proper closed subsets of $X$ such that $F_i\cap F_j\in\mathcal{C}$
for all $i\neq j$. Let $F=\bigcup_{i\neq j}F_i\cap F_j$. Shrinking
$F_i$, $i\geq 0$, to smaller closed subsets and re-indexing these
sets, if necessary, we can assume that there exist $n\neq m$ and two
closed disjoint subsets $X_0$ and $X_1$ of $X$ both having non-empty
interiors in $X$ with $X_0\subset F_n\setminus\bigcup_{i\neq n}F_i$
and $X_1\subset F_m\setminus\bigcup_{i\neq m}F_i$. This can be done
using arguments similar to the Baire theorem. Then $X_0\cup X_1$ is
disjoint from the set $F$. Since $X$ is a Mazurkiewicz manifold with
respect to $\mathcal{C}$, there exists a continuum $C\subset
X\setminus F$ joining $X_0$ and $X_1$. This implies that $C$ is
covered by the family $F_i\cap C$ of its disjoint closed subsets.
Hence, according to the Sierpi\'{n}ski theorem~\cite[p. 440]{re:89},
$C=F_i\cap C$ for some $i$, which contradicts the conditions on
$X_0$ and $X_1$.
\end{proof}

The following lemma is a variation on the countable sum theorem and
is required in the proof of Proposition \ref{extension1}.

\begin{lem}\label{extensionfsigma}
Let $X$ be a compact space, $A$ be a closed subspace of $X$, and $Y=
\bigcup_{n\geq 1}Y_n$ with all $Y_n$ being closed in $X$. Let $L$ be
a $CW$-complex such that $L\in AE(Y_n)$, $n\geq 1$. Then any map
$f\colon A\to L$ is extendable over some open neighborhood of $A\cup
Y$ in $X$.
\end{lem}

\begin{proof}
We may assume that $A\subset Y_{n}\subset Y_{n+1}$ for $n\geq 1$.
Let $f_0 = f$ and construct by induction a sequence of maps
$f_n\colon\cl(U_n) \to L$, where $U_n$ is an open neighborhood of $
Y_n$ in $X$. Suppose $f_n$ has been already constructed. Since $L\in
AE(Y_{n+1})$, we can extend $f_n$ to a map $\overline{f}_n
\colon\cl(U_n)\cup Y_{n+1}\to L$. Using that $L$ is an absolute
neighborhood extensor for compact spaces, we extend $\overline{f}_n$
to a map $f_{n+1}\colon\cl(U_{n+1}) \to L$, where $U_{n+1}$ is an
open neighborhood of $\cl(U _n) \cup Y_{n+1}$. The sequence $\{
f_n\}_{n\ge 0}$ gives rise to a map
$$\overline{f}\colon U = \bigcup _{n\ge 0} U_n \to L$$
extending $f$.
\end{proof}

Propositions \ref{extension1} and \ref{extension2}  below are among
the main technical tools. As we noted in the Introduction,
particular cases of Proposition~\ref{extension2} were established by
various authors.

\begin{pro}\label{extension1}
Let $L$ be a $CW$-complex and $X$ be a compact space. Let
$\{F_i\}_{i\ge 0}$ be a family of closed subsets of $X$ such that
$L\in AE \big(F_i\times\I\big)$ for all $i$ and $F = \bigcup _{i\ge
0} F_i$ cuts $X$ between two closed subsets $X_0$ and $X_1$ both
having non-empty interiors. Let $A\subset X$ be a closed set and
$f\colon A\to L$ be a map extendable over $A\cup Y$ for every proper
closed subset $Y$ of $X$. Then $f$ is extendable over $X$.
\end{pro}

\begin{proof}
We may assume that $X_0 =\cl(U_0)$ and $X_1 =\cl(U_1)$, where $U_0$
and $U_1$ are non-empty open subsets of $X$. Let $Y_0 = X\setminus
U_1$ and $Y_1 = X\setminus U_0$. Then both $Y_0$ and $Y_1$ are
proper closed subsets of $X$. Therefore there exist two maps
$f_0\colon Y_0\cup A \to L$ and $f_1\colon Y_1 \cup A \to L$ both
extending $f$.

Consider the map $G\colon (Y_0 \times \{0\}) \cup
(Y_1\times\{1\})\cup A\times\I\to L$ defined as follows:

$$G(x,t) = \left\{
        \begin{array}{ll}
          f_0(x), & \hbox{if $x\in Y_0$ and $t=0$;} \\
         f_1 (x), & \hbox{if $x\in Y_1$ and $t=1$;} \\
          f(x), & \hbox{if $x\in A$.}
        \end{array}
      \right.
$$

According to Lemma~\ref{extensionfsigma}, the map $G$ can be
extended to a map $H\colon W\to L$, where $W$ is an open
neighborhood of $(Y_0 \times \{0\}) \cup (Y_1\times\{1\})\cup
(A\times\I)\cup (F\times\I)$ in $X\times\I$. Since $\I$ is compact,
there is an open set $V\subset X$ containing $F$ such that
$V\times\I \subset W$ and $V\cap (X_0\cup X_1)=\varnothing$.

If $X\setminus V$ were connected between $X_0$ and $X_1$, then there
would be a continuum $C\subset X\setminus V$ such that $C\cap X_k
\ne \varnothing$, $k=0,1$ (see, e.g., \cite[ \S 47.II, Theorem
3]{Ku1}),  contradicting the fact that $F$ cuts $X$ between $X_0$
and $X_1$. Therefore, the set $V$ contains a closed partition $P$
between $X_0$ and $X_1$ in $X$. Thus $X = X' _0 \cup X' _1$, where
$X'_0$ and $X'_1$ are closed subsets of $X$ such that $X'_0\cap X'_1
= P$ and $X_k \subset X'_k$, $k=0,1$. According to the definition of
$Y_0$ and $Y_1$, we have $X'_k \subset Y_k$, $k=0,1$.

Let $f'_k = f_k |_{A\cup P}$, $k=0,1$. Note that the map $H|_{(A\cup
P)\times\I}$ is a homotopy between $f'_0$ and $f'_1$. Then, by the
Homotopy Extension Theorem, there exists a map from $X$ into $L$
extending $f$.
\end{proof}

\begin{pro}\label{extension2}
Let $L$ be a $CW$ complex and $X$ be a compact space admitting a
cover $\{F_i\}_{i\ge 0}$ by closed subsets $F_i\subset X$ such that
$L\in AE\big((F_i\cap F_j)\times\I\big)$ for all $i\ne j$. Let
$A\subset X$ be a closed set and $f\colon A\to L$ a map extendable
over $A\cup F_i$ for every $i$. Then $f$ is extendable over $X$.
\end{pro}

\begin{proof}
Suppose the opposite. Let $\mathcal{A}$ be the family of all closed
subsets $Y$ of $X$ containing $A$ such that $f$ is not extendable
over $Y$. Note that $\mathcal{A}$ is partially ordered by inclusion
and $X\in\mathcal{A}$. We show that $\mathcal A$ satisfies the
Zorn's lemma. Indeed, suppose $\{Y_\alpha:\alpha\in\Lambda\}$ is a
decreasing net of sets from $\mathcal A$ and
$Y=\bigcap\{Y_\alpha:\alpha\in\Lambda\}$ is not in $\mathcal A$. If
there exists a map $\overline{f}\colon Y\to L$ extending $f$, then
$\overline{f}$ can be extended to a map $g\colon U\to L$ with $U$
being an open neighborhood of $Y$ in $X$. Due to the compactness,
$U$ contains $Y_\alpha$ for some $\alpha$, which is a contradiction.

Let $M$ be a minimal element of $\mathcal A$. Let $C_i = M\cap F_i$.
Since $f$ is extendable over each $A\cup F_i$ but not extendable
over $M$, all $A\cup C_i$, $i\geq 0$, are proper subsets of $M$.
Using this fact and the Baire theorem, we can assume that there
exist open sets $U_0$ and $U_1$ in $M$ such that

\begin{align}
& \cl(U_0)\subset C_0\setminus A, \quad \cl(U_1)\subset C_1\setminus A, \quad  \cl(U_0)\cap\cl(U_1) = \varnothing \notag\\
& \text{and}\notag\\
& U_0 \cap C_0 \cap C_1 =\varnothing, \quad U_1 \cap C_0\cap C_1
=\varnothing.\notag
\end{align}

Denote
$$B_0 = C_0, \quad B_1 = C_1, \quad B_i = C_i \setminus (U_0\cup
U_1)\quad\text{for $i \ge 2$}$$ and let  $B = \cup _{i\ne j}(B_i\cap
B_j)$. Shrinking $U_0$ and $U_1$, if necessary, we may also assume
that
$$\cl(U_0)\cap \cl(B)
=\varnothing \quad\text{and}\quad  \cl(U_1)\cap \cl(B)
=\varnothing.$$

We claim that $B$ cuts $M$ between $\cl(U_0)$ and $\cl(U_1)$.
Indeed, suppose not. Then there exists a continuum $C\subset
M\setminus B$ such that $C\cap\cl(U_k) \ne \varnothing$, $k=0,1$.
Note that $\{B_i\cap C\}_{i\ge 0}$ is a cover of $C$ by closed
disjoint proper sets. Hence, by the Sierpi\'{n}ski theorem \cite[p.
440]{re:89}, $C \subset B_i$ for some $i$, which contradicts the
choice of $\cl(U_0)$ and $\cl(U_1)$.

Therefore, due to the minimality of $M$, we can apply Proposition
\ref{extension1} to $M$, the collection $B_i\cap B_j$, $i,j\ge 0$,
$i\ne j$, and the sets $\cl(U_0)$ and $\cl(U_1)$ to obtain a map
$\overline{f}\colon M \to L$ extending $f$. This contradicts
$M\in\mathcal{A}$.
\end{proof}

The following technical lemma will help us to work with Mazurkiewicz
manifolds.

\begin{lem}\label{mazur1}
Let $X$ be a compact space, $X_0$ and $X_1$ be two closed disjoint
subsets of $X$ with non-empty interiors, and $S$ be a subset of $X$.
Suppose that for any continuum $C$ with $C\cap X_0 \ne
\varnothing\ne C\cap X_1$ we have $C\cap S \ne\varnothing$. Then
there exist open non-empty sets $U_k$ and $V_k$ with
$V_k\subset\cl(V_k)\subset U_k \subset X_k$, $k=0,1$, such that for
any continuum $C$ with $C\cap\cl(V_0) \ne \varnothing\ne
C\cap\cl(V_1)$ we have $C\cap (S\setminus (U_0\cup U_1))\ne
\varnothing$.
\end{lem}

\begin{proof} Since $X_0$ and $X_1$ have non-empty interiors, we can find
open non-empty sets $U_k$ and $V_k$ such that $\cl(V_k)\subset
U_k\subset X_k$, $k=0,1$. Consider a continuum $C$ such that
$C\cap\cl(V_0) \ne \varnothing\ne C\cap\cl(V_1)$. Note that
$C\cap\bd(U_k) \ne\varnothing$, $k=0,1$. Since $C$ is a continuum,
there exists a  component $C'$ of the compact space $C\setminus
(U_0\cup U_1)$ such that $C'\cap\bd(U_k) \ne\varnothing$, $k=0,1$.
Then $C'$ is a continuum joining $X_0$ and $X_1$ and therefore
$C'\cap S\ne\varnothing$. Since $C'\subset C\setminus (U_0\cup
U_1)$, we have $C\cap (S\setminus (U_0\cup U_1))\ne \varnothing$, as
required.
\end{proof}

Now we are ready to prove our first main result.

\begin{thm}\label{main1}
Every compact space $X$ with $D_{\mathcal{K}}(X)=n\geq 1$ contains a
closed subset $M$ such that $D_{\mathcal{K}}(M)=n$ and $M$ is both a
$V^n$-continuum and a Mazurkiewicz manifold (and hence a strong
Cantor $n$-manifold) with respect to $\mathcal
D_{\mathcal{K}}^{n-2}$.
\end{thm}

\begin{proof}
Since $D_{\mathcal K}(X)=n$, we have $K_n\in AE(X)$ but
$K_{n-1}\not\in AE(X)$. Therefore there exists a closed subset
$A\subset X$ and a map $f\colon A\to K_{n-1}$ which cannot be
extended to a map from $X$ into $K_{n-1}$. Consider the family
$\mathcal{B}$ of all closed sets $B\subset X$ such that there is no
map from $A\cup B$ to $K_{n-1}$ extending $f$. Obviously,
$X\in\mathcal{B}$. As in the proof of Proposition \ref{extension2},
one verifies that $\mathcal{B}$ is partially ordered by inclusion
and satisfies the condition of the Zorn's lemma. Let $M$ be a
minimal element of $\mathcal{B}$. Then, $D_{\mathcal K}(M)\leq
D_{\mathcal K}(X)= n$. Since the map $f|_{A\cap M}$ cannot be
extended to a map from $M$ into $K_{n-1}$, $D_{\mathcal K}(M)>n-1$.
Thus, $D_{\mathcal K}(M)=n$.

Suppose $M$ is not a $V^n$-continuum with respect to $D_{\mathcal
K}^{n-2}$. Then, without loss of generality, we can assume that
there exist two disjoint sets $X_0 =\cl(U_0)\subset M$ and $X_1 =
\cl(U_1) \subset M$, where $U_0$ and $U_1$ are open non-empty
subsets of $M$,  with the following property:

\begin{itemize}\item
for any open cover $\omega$ of $M$ there exists a partition
$P_\omega$ in $M$ between $X_0$ and $X_1$ such that $P_\omega$
admits an $\omega$-map into a compact space $Y_\omega$ of dimension
$D_{\mathcal K}(Y_\omega)\leq n-2$.
\end{itemize}

Since both $M_0= M\setminus U_1$ and $M_1=M\setminus U_0$ are proper
closed subsets of $M$  and $M$ is a minimal element of
$\mathcal{B}$, there exist maps $f_i\colon M_i\to K_{n-1}$ extending
$f| _{A\cap M_i}$, $i=0,1$. Let
$$Z=(M_0\times\{0\})\cup (M_1\times \{1\})\cup\big((M\cap A)
\times \I\big).$$

Define a map $F\colon Z\to K_{n-1}$ by

$$F(x,t) = \left\{
        \begin{array}{ll}
          f_0(x), & \hbox{if $x\in M_0$ and $t=0$;} \\
         f_1 (x), & \hbox{if $x\in M_1$ and $t=1$;} \\
          f(x), & \hbox{if $x\in A\cap M$.}
        \end{array}
      \right.
$$

Let $\gamma$ be an open cover of $K_{n-1}$ such that any two
$\gamma$-close maps to $K_{n-1}$ are homotopic.

Next claim follows easily from the fact that $K_{n-1}$, as a
metrizable $ANR$, is a neighborhood retract of a locally convex
space (see \cite[Lemma 8.1]{bv} for a similar proof).

\begin{claim} There exists an open cover $\nu$ of $Z$ satisfying
the following condition: for any closed $B\subset Z$ and any
$\nu$-map $\varphi\colon B\to Y$ into a paracompact space $Y$, there
exists a map $g\colon \varphi (B)\to K_{n-1}$ such that $F|_B$ and
$g\circ\varphi$ are $\gamma$-close in $K_{n-1}$.
\end{claim}

Let $\omega$ be an open cover of $M$ such that each set $(W\times
\{t\})\cap Z$, $W\in\omega$ and $t\in\I$, is contained in some
element of $\nu$. There exists a partition $P_\omega$ in $M$ between
$X_0$ and $X_1$ admitting an $\omega$-map $\varphi_\omega\colon
P_\omega\to Y_\omega$ into a compact space $Y_\omega$ of dimension
$D_{\mathcal K}(Y_\omega)\leq n-2$.

Let
$$B=P_{\omega}\times\{0,1\}\cup\big((P_{\omega}\cap A) \times
\I\big)$$
 and $\varphi\colon B\to Y = Y_{\omega}\times\I$ be defined
as
$$\varphi (x,t) = (\varphi_{\omega} (x),t)\quad\text{for all $(x,t)\in  B$}.$$
 Note that $\varphi$ is a $\nu$-map. Applying the above claim we
obtain a map $g\colon \varphi (B)\to K_{n-1}$ such that $F|_B$ and
$g\circ\varphi$ are $\gamma$-close in $K_{n-1}$.  The map
$$\Phi\colon P_{\omega}\times\I\to Y_{\omega}\times\I,\quad \Phi (x,t) =
(\varphi _{\omega}(x),t),$$
 is an extension of $\varphi$. Since
$D_{\mathcal K}(Y_\omega\times\I)\leq n-1$, the map $g$ can be
extended to a map $G\colon Y_\omega\times\I\to K_{n-1}$. Note that
$F|_B$ and $(G\circ\Phi)|_B = g\circ\varphi$ are $\gamma$-close, and
therefore homotopic by the choice of $\gamma$. The Homotopy
Extension Theorem implies the existence of a map $H\colon
P_{\omega}\times\I\to K_{n-1}$ extending $F|_B$. Note that $H$ is a
homotopy between $f_0|_{P_{\omega}}$ and $f_1|_{P_{\omega}}$ such
that $H(x,t) = f(x)$ for all $x\in P_{\omega}\cap A$. Since
$P_{\omega}$ is a partition between $X_0$ and $X_1$, there exist two
closed subsets $M'_0$ and $M'_1$ of $M$ such that $X_i\subset
M'_i\subset M_i$, $i=0,1$, $M'_0\cup M'_1 = M$ and $M'_0\cap M'_1 =
P_{\omega}$. Applying the Homotopy Extension Theorem to the space
$M'_1$, its closed subset $P = P_{\omega}\cup (A\cap M'_1)$, and the
maps $f_0$ and $f_1$, we get a map $f'_0\colon M'_1\to K_{n-1}$
extending $f_1|_{P}$ over $M'_1$. By pasting $f_0$ and $f'_0$ we
finally obtain an extension of $f|_{M\cap A}$ over $M$. This yields
a contradiction with $M\in\mathcal{B}$. Thus, $M$ is a
$V^n$-continuum with respect to $D_{\mathcal K}^{n-2}$.

Now we show that $M$ is a Mazurkiewicz manifold with respect to
$\mathcal D_{\mathcal K}^{n-2}$. Assuming the opposite and applying
Lemma \ref{mazur1}, we find  closed subsets $F_i$  of $M$, $i\geq
0$, such that $F=\bigcup _{i\ge 0} F_i$ cuts $M$ between two closed
disjoint subsets of $M$ with non-empty interiors and  $D_{\mathcal
K}(F)\leq n-2$.

Note that $K_{n-2}\in AE(F_i)$ for each $i$. So, according to the
definition of a stratum, $K_{n-1}\in AE(F_i\times\I)$. Moreover,
since $M$ is a minimal element of $\mathcal{B}$, the map $f|(A\cap
M)$ can be extended to a map from $(A\cap M)\cup Y$ into $K_{n-1}$
for any proper closed subset $Y$ of $M$. Then, by
Proposition~\ref{extension1}, there exists a map $g\colon M\to
K_{n-1}$ extending $f|(A\cap M)$, which contradicts
$M\in\mathcal{B}$.

\end{proof}


\section{Infinite-dimensional Mazurkiewicz manifolds}

In this section we consider Mazurkiewicz manifolds with respect to
classes $\mathcal D_{\mathcal{K}}^{<\infty}$, $\mathcal{WID}$ and
$\mathbf C$ (of strongly countable $\mathcal
D_{\mathcal{K}}$-dimensional spaces, weakly infinite-dimensional
spaces and $C$-spaces, respectively).

\begin{thm}\label{inf}
If a compact space $X$ has dimension $D_{\mathcal{K}}(X)=\infty$,
then either $X$ contains closed subsets of arbitrary large finite
dimensions  $D_{\mathcal{K}}$ or $X$ contains a compact Mazurkiewicz
manifold with respect to  the class $\mathcal
D_{\mathcal{K}}^{<\infty}$.
\end{thm}

\begin{proof}
We have $K_n\notin AE(X)$ for all $n\ge0$. Suppose there exists
$n_0\in\mathbb N$ such that $X$ contains no closed subset of
dimension $D_{\mathcal{K}}\ge n_0$. We  follow the idea from the
proof of  Theorem~\ref{main1}. First, choose a closed subset
$A\subset X$ and a map $f\colon A\to
 K_{n_0}$ which cannot be extended over $X$. Then, there exists $M$  minimal
in the family $\mathcal B$ of all closed subsets $B\subset X$ for
which there is no extension of $f$ over $A\cup B$. It follows that
$D_{\mathcal{K}}(M)\ge n_0+1$, hence  $D_{\mathcal{K}}(M)=\infty$.

Suppose  $M$ is not a Mazurkiewicz  manifold with respect to the
class $\mathcal D_{\mathcal{K}}^{<\infty}$. Then, by
Lemma~\ref{mazur1}, there exist closed subsets $F_i\subset M$ such
that $F=\bigcup_{i\ge 1}F_i$ cuts $M$ between two closed, disjoint
subsets of $M$ with non-empty interiors and
$D_{\mathcal{K}}(F)=n<\infty$ for some $n<n_0$. It follows that
$K_n\in AE(F_i)$, so $\displaystyle K_{n+1}\in AE(F_i\times\mathbb
I)$ for each $i$. Since $n+1\leq n_0$, $\displaystyle K_{n_0}\in
AE(F_i\times\mathbb I)$, $i\geq 1$. The minimality of $M$ implies
that the map $f|(A\cap M):A\cap M \to\displaystyle K_{n_0}$ extends
over $(A\cap M)\cup Y$ for any proper closed subset $Y\subset M$.
Now, by Proposition~\ref{extension1}, there exists an extension of
$f|(A\cap M)$ over $M$, a contradiction  with $M\in \mathcal B$.
 \end{proof}

 Recall that a set-valued map $\Phi\colon X\to Y$ is lower
semi-continuous (resp., upper semi-continuous) if the set $\{x\in
X:\Phi(x)\cap U\neq\varnothing\}$ (resp., $\{x\in X:\Phi(x)\subset
U\}$) is open in $X$ for every open $U\subset Y$. We say that $\Phi$
is continuous provided it is both lower semi-continuous and upper
semi-continuous. Recall also that a closed subset
$F\subset\I^\infty$ is said to be a $Z$-set in $\I^\infty$ if for
every compact space $X$ the set $\{g\in C(X,\I^\infty): g(X)\cap
F=\varnothing\}$ is dense in $C(X,\I^\infty)$ in the compact-open
topology.

\begin{pro}\label{pro4}
A compact space $X$ does not have property $C$ if and only if there
exists a continuous set-valued map $\Phi\colon X\to\I^\infty$
satisfying the following conditions: each $\Phi(x)$ is a $Z$-set in
$\I^\infty$ and for any single-valued map $g\colon X\to\I^\infty$ we
have $g(x)\in\Phi(x)$ for some $x\in X$.
\end{pro}

\begin{proof}
This proposition is a direct consequence of the following result of
Uspenskij \cite[Theorem 1.4]{vu} that characterizes compact
$C$-spaces: a compact space  has the property $C$ if and only if for
every continuous $\Phi\colon X\to\I^\infty$ with each $\Phi(x)$
being a $Z$-set in $\I^\infty$ there exists a single-valued map
$g\colon X\to\I^\infty$ such that $g(x)\not\in\Phi(x)$ for all $x\in
X$.
\end{proof}

\begin{lem}\label{lem3}
Let $X$ be a compact space and $\Phi\colon X\to\I^\infty$ a
continuous set-valued map with each $\Phi(x)$ being a $Z$-set in
$\I^\infty$. Suppose $A\subset X$ is closed and $F=\bigcup_{i\geq 1}
F_i$ such that all $F_i$ are closed $C$-subspaces of $X$. Then any
map $f\colon A\to\I^\infty$ with $f(x)\not\in\Phi(x)$, $x\in A$, can
be extended to a map $g\colon W\to\I^\infty$, where $W$ is a
neighborhood of $A\cup F$, such that $g(x)\not\in\Phi(x)$ for any
$x\in W$.
\end{lem}

\begin{proof}
Consider the sets $$C(f)=\{h\in C(X,\I^{\infty}):h|A=f\}$$
 and
$$C_{i}(f)=\{h\in C(f):\text{$h(x)\not\in\Phi(x)$ for
all $x\in F_{i}$}\}, \quad i\geq 1.$$ Here, $C(X,\I^\infty)$ is the
space of all continuous maps from $X$ into $\I^\infty$ equipped with
the metric $d(g_1,g_2)=\max\{\rho\big(g_1(x),g_2(x)\big):x\in X\}$,
where $\rho$ is the standard convex metric on $\I^\infty$.

We claim that each $C_{i}(f)$ is open and dense in $C(f)$. Indeed,
let $h\in C_{i}(f)$ and observe that
$\epsilon=\min\{\rho\big(h(x),\Phi(x)\big):x\in F_{i}\}$ is positive
because $\Phi$ is continuous. Then, any map in $C(f)$ which is
$\epsilon$-close to $h$ is contained in $C_{i}(f)$. Thus $C_{i}(f)$
is open in $C(f)$.

To prove $C_{i}(f)$ is dense in $C(f)$, fix $h\in C(f)$ and
$\epsilon=2\eta>0$, and consider the set-valued map
$$\phi\colon
F_{i}\to\I^{\infty},\quad
\phi(x)=\begin{cases} f(x)\quad&\text{for $x\in A\cap F_i$,}\\
                        \overline{B}\big(h(x),\eta\big) \quad&\text{for $x\in F_i\setminus
A$},
\end{cases}$$
where $\overline{B}\big(h(x),\eta\big)$ is the closed ball in
$(\I^{\infty},\rho)$ with radius $\eta$ and center $h(x)$. This is a
lower semi-continuous convex-valued map. Since all
$\overline{B}\big(h(x),\eta\big)$ are convex and $\Phi(x)$ are
$Z$-sets in $\I^{\infty}$, it is easily seen that
$\overline{B}\big(h(x),\eta\big)\cap\Phi(x)$ is a $Z$-set in
$\overline{B}\big(h(x),\eta\big)$, $x\in F_{i}\setminus A$. Since
$F_{i}$ is a $C$-space, by \cite[Theorem 1.1]{gv:02}, $\phi$ admits
a continuous selection
$$h_1\colon F_{i}\to\I^{\infty}\quad\text{with}\quad
h_1(x)\not\in\Phi(x) \quad\text{for all}\quad x\in F_{i}.$$
 Now, define
$$h_2\colon
A\cup F_i\to\I^\infty \quad\text{by}\quad h_2|A=f
\quad\text{and}\quad h_2|F_i=h_1.$$ Finally, extend $h_2$ to a map
$h_3\in C(X,\I^{\infty})$ in such a way that $h_3$ is $\eta$-close
to $h$. According to the convex-valued selection theorem of Michael
\cite{m}, the map $h_3$ can be obtained as a selection of the
convex-valued lower semi-continuous map
$$\varphi\colon X\to\I^{\infty},\quad \varphi(x)=\begin{cases} h_2(x)\quad&\text{if $x\in A\cup
F_{i}$,}\\
\overline{B}\big(h(x),\eta\big)\quad&\text{otherwise}.
\end{cases}$$
Obviously, $h_3\in C_{i}(f)$ and it is $\epsilon$-close to $h$.

Since $C(f)$ is complete (as a closed subset of $C(X,\I^\infty)$),
by the Baire theorem, there exists a map $g\in\bigcap_{i\geq
1}C_{i}(f)$. Then $g(x)\not\in\Phi(x)$ for all $x\in F\cup A$ and
$g|A=f$. Moreover, by the continuity of $\Phi$, one can show that
every point $x\in F\cup A$ has a neighborhood $O(x)$ in $X$ with
$g(y)\not\in\Phi(y)$ for all $y\in O(x)$. Then $W=\bigcup_{x\in
F\cup A}O(x)$ is a neighborhood of $A\cup F$ such that
$g(x)\not\in\Phi(x)$ for $x\in W$.
\end{proof}

\begin{thm}\label{C}
Every compact space $X$ which is not a $C$-space contains a compact
Mazurkiewicz manifold with respect to the class $\mathbf C$.
\end{thm}

\begin{proof}
Let $\Phi\colon X\to\I^\infty$ be a continuous set-valued map
satisfying Proposition~\ref{pro4}. Consider the family
$\mathcal{B}_{\Phi}$ of all closed subsets $B\subset X$ such that
for every map $g\colon B\to\I^\infty$ there exists a point $x\in B$
with $g(x)\in\Phi(x)$. Let us show that $\mathcal{B}_{\Phi}$ has a
minimal element. Indeed, if $\{B_\alpha:\alpha\in\Lambda\}$ is a
decreasing net of sets from $\mathcal{B}_{\Phi}$ and
$B_0=\bigcap\{B_\alpha:\alpha\in\Lambda\}$, then every $g\colon
B_0\to\I^\infty$ can be extended to a map $\overline{g}\colon
X\to\I^\infty$. For every $\alpha\in\Lambda$ choose $x_\alpha\in
B_\alpha$ such that $\overline{g}(x_\alpha)\in\Phi(x_\alpha)$ and
let $x_0$ be a limit point of a subnet of $\{x_\alpha\}$. Obviously,
$x_0\in B_0$ and since both $\Phi$ and $\overline{g}$ are
continuous, $g(x_0)\in\Phi(x_0)$. Thus, by the Zorn lemma,
$\mathcal{B}_{\Phi}$ has a minimal element $M$. Since
$M\in\mathcal{B}_\Phi$, Proposition~\ref{pro4} yields that $M$ is
not a $C$-space.


We will show that $M$ is a Mazurkiewicz manifold with respect to the
class $\mathbf C$. Suppose not. Then, by Lemma \ref{mazur1}, there
exist closed subsets $F_i$, $i=0,1,2,\dots$,  of $M$ such that
$F=\bigcup _{i\ge 0} F_i$ cuts $M$ between two closed disjoint
subsets $X_0 =\cl(U_0)$ and $X_1 =\cl(U_1)$, where $U_0$ and $U_1$
are non-empty open subsets of $M$, and $F$ is a $C$-space.

Let $Y_0 = M\setminus U_1$ and $Y_1 = M\setminus U_0$. Then both
$Y_0$ and $Y_1$ are proper closed subsets of $M$. Therefore there
exist two maps $g_i\colon Y_i\to \I^\infty$ such that $g_i(x) \notin
\Phi (x)$ for all $x\in Y_i$, $i=0,1$. Consider the map $g\colon
(Y_0 \times \{0\}) \cup (Y_1\times\{1\})\to\I^\infty$ defined as
follows:

$$g(x,t) = \left\{
        \begin{array}{ll}
          g_0(x), & \hbox{if $x\in Y_0$ and $t=0$;} \\
         g_1 (x), & \hbox{if $x\in Y_1$ and $t=1$.} \\
                 \end{array}
      \right.
$$

\noindent Applying Lemma~\ref{lem3} to the closed subset $A = (Y_0
\times \{0\}) \cup (Y_1\times\{1\})$ of $M\times\I$ and to
$F\times\I$ (which is a $C$-space), we obtain an extension $G\colon
W\to \I^{\infty}$ of $g$ over some open neighborhood $W$ of $A\cup
(F\times \I)$ in $M\times\I$, such that $G(x,t)\notin \Phi (x)$ for
all $(x,t)\in W$. Due to the compactness of $\I$, we can find an
open subset $V$ of $M$ containing $F$ such that $V\times\I\subset W$
and $V\cap (X_0\cup X_1)= \varnothing$. As in the proof of
Proposition~\ref{extension1} we conclude that $V$ is an open
partition between $X_0$ and $X_1$ in $M$. Then $M\setminus V =
M_0\cup M_1$, where $M_i$ are disjoint closed subsets of $M$ and
$X_i\subset M_i\subset Y_i$, $i=0,1$. Let $\theta\colon M\to\I$ be a
function such that $\theta (M_i)=i$, $i=0,1$. Then the map $f(x) =
G(x,\theta(x))$ is well-defined for all $x\in M$ and
$f(x)\notin\Phi(x)$ for any $x\in M$. The last condition contradicts
$M\in\mathcal{B}_\Phi$.

Thus, $M$ is a Mazurkiewicz manifold with respect to the spaces
having property $C$.
\end{proof}

The next theorem is an analogue of Theorem~\ref{C} for strongly
infinite-dimensional spaces. We say that a (single-valued) map
$f\colon X\to\I^\infty$ is {\it universal} \cite{ho1} if for any map
$g\colon X\to\I^\infty$ there exists a point $x\in X$ with
$g(x)=f(x)$.

\begin{pro}
A compact space $X$ is strongly infinite-dimensional if and only if
there exists a universal map $f\colon X\to\I^\infty$.
\end{pro}

\begin{proof}
By \cite{ap}, a compact space $X$ is strongly infinite-dimensional
if and only if there exists an essential map $f\colon
X\to\I^\infty$. Recall that a map $f\colon X\to\I^\infty$ is
essential if for every $n$ the composition $\pi_n\circ f$ is
essential, i.e. there is no map $g\colon X\to\mathbb{S}^{n-1}$ with
$g|(\pi_n\circ f)^{-1}(\mathbb{S}^{n-1})=(\pi_n\circ f)|(\pi_n\circ
f)^{-1}(\mathbb{S}^{n-1})$. Here, $\pi_n\colon\I^\infty\to\I^n$ is
the projection onto $\I^n$ and $\mathbb{S}^{n-1}$ is the boundary of
$\I^n$. On the other hand, a map $f\colon X\to\I^\infty$ is
essential if and only if $f$ is universal (this fact was established
in \cite{gr} for metrizable compact spaces, but the proof works for
arbitrary compact spaces).
\end{proof}

\begin{thm}\label{wid}
Every strongly infinite-dimensional metrizable compact space $X$
contains a Mazurkiewicz manifold with respect to the class
$\mathcal{WID}$.
\end{thm}

\begin{proof}
We fix a (single-valued) universal map $\Phi\colon X\to\I^{\infty}$.
Observe that the values of $\Phi$, being points, are $Z$-sets in
$\I^\infty$. Since $X$ contains a strongly infinite-dimensional
closed set $Y$ such that every subset of $Y$ is either
$0$-dimensional or strongly infinite-dimensional (see \cite{lr} or
\cite{ml}), we can assume that every subset of $X$ is
$0$-dimensional provided it is not strongly infinite-dimensional.
Then, as in the proof of Theorem~\ref{C}, we can obtain a closed
strongly infinite-dimensional set $M\subset X$ such that the map
$\Phi|M\colon M\to\I^\infty$ is universal, but $\Phi|H$ is not
universal for any closed proper subset $H$ of $M$. Following the
ideas  from the proof of Theorem~\ref{C}, we can show that $M$ is a
Mazurkiewicz manifold with respect to the class $\mathcal{WID}$.
Indeed, if $\{F_i\}_{i\geq 1}$ is a sequence of closed subsets of
$M$ with $F_i$ being weakly infinite-dimensional, then $F_i$ should
be $0$-dimensional. Note that every $0$-dimensional compact space is
a $C$-space, so we can apply the arguments from the proof of
Theorem~\ref{C}.
\end{proof}

\section{Applications to homogeneous continua}

All spaces in this section are metrizable and the dimension of a
space $X$ means any dimension $D_{\mathcal K}(X)$ if not stated
otherwise.

\begin{re}\label{r}
Recall that a connected, locally compact metrizable space $X$ is
second-countable. Thus, by the Countable Sum Theorem, if $X$
contains a closed $n$-dimensional subset, then $X$ contains compact
$n$-dimensional subsets of arbitrary small   diameters.
\end{re}

A topological group $H$ acts transitively on a space $X$ if the
action $H\times X\to X$ is continuous and for each two points
$x,y\in X$ there is $h\in H$ such that $h(x)=y$. We denote by $H(X)$
the group of homeomorphisms of a space $X$ onto itself with a
compact-open topology. A space $X$ is homogeneous if $H(X)$ acts
transitively on $X$, i.e. for  each two points $x,y\in X$ there
exists $h\in H(X)$ such that $h(x)=y$;
 $X$ is called locally homogeneous if for  each  $x,y\in X$ there exist
 neighborhoods $U$ and $V$ of $x$ and $y$, respectively,
and a homeomorphism $h:U\to V$ such that $h(x)=y$.

\begin{thm}[Effros' Theorem~\cite{Ef}]\label{Effros} If  $H(X)$ acts transitively on
a closed subset $Y$ of a compact space $(X,d)$, then for every
$\epsilon>0$ there exists $\delta>0$ such that if $x,y\in Y$ and
$d(x,y)<\delta$, then there exists $h\in H(X)$ such that $h(x)=y$
and $d(h(z),z)<\epsilon$ for every $z\in X$.
\end{thm}

A homeomorphism $h$ in the above theorem will be called an
$\epsilon$-ho\-meo\-morphism.

\

The following simple observation explains a role of the class
$\mathcal D_{\mathcal K}^{<\infty}$ for infinite-dimensional
homogeneous continua.

\begin{pro}\label{infty}
A homogeneous continuum is infinite-dimensional if and only if it is
not strongly countable dimensional.
\end{pro}

\begin{proof}
Suppose $X$ is a  homogeneous continuum and $X=\bigcup_{i=1}^\infty
F_i$, where $F_i$ is a closed finite-dimensional closed subset of
$X$ for each $i$. There exists $k$ such that $\inte
F_k\ne\emptyset$, by the Baire theorem. By the homogeneity, finitely
many homeomorphic copies of $F_k$ covers continuum $X$,  so it is
finite-dimensional. The converse implication is obvious.
\end{proof}

\begin{thm}\label{TC}
Each  homogeneous continuum $X\notin\mathcal{C}$ is a  Cantor
manifold with respect to class $\mathcal{C}$ where $\mathcal{C}$ is
any of the  following four classes: $\mathcal D_{\mathcal
K}^{<\infty}$, $\mathcal{WID}$, $\mathbf C$, $\mathcal D_{\mathcal
K}^{n-2}$
 $($in the latter case  we additionally assume
$D_{\mathcal{K}}(X)=n$$)$.
\end{thm}

\begin{proof} Theorem~\ref{TC} was proved in~\cite{kru} for the covering
dimension and weak infinite dimension. The proof was based on the
classical Cantor Manifold Theorem  that any compact $n$-dimensional
space contains a  Cantor $n$-manifold (see~\cite{re:95}), the
corresponding Tumarkin's result for infinite-dimensional compacta
and Skljarenko's theorem for the case of strongly
infinite-dimensional compacta (both  mentioned in the Introduction).
Due to Theorems~\ref{main1}, \ref{inf}, \ref{C},  the same idea
applies. In the case of class $\mathcal{C}=\mathcal D_{\mathcal
K}^{<\infty}$, however, we have to consider an extra   situation
when there is no Cantor manifold with respect to $\mathcal
D_{\mathcal K}^{<\infty}$ in $X$ but $X$ contains closed subsets of
arbitrary large finite dimension  (see Theorem~\ref{inf}). In
particular, $X$ is not a Cantor manifold with respect to $\mathcal
D_{\mathcal K}^{<\infty}$ which means that there is a closed set
$F=\bigcup_{n=1}^\infty F_n$ which separates $X$,  where $F_n$ is a
finite-dimensional closed set. We can assume that $X\setminus
F=U\cup V$, where $U, V$ are non-empty disjoint open subsets of $X$
and $F=\bd U =\bd V$. By the Baire theorem, one can find $n_0$ such
that $\inte_F (F_{n_0})\ne\emptyset$. Since each  finite-dimensional
nondegenerate compactum  contains arbitrary small   compacta of the
same finite dimension (Remark~\ref{r}), there are arbitrary small
Cantor manifolds in $X$ with respect to this finite dimension
(Theorem~\ref{main1}). Let $D_{\mathcal{K}}(F_{n_0})=n$ and pick up
a point $x\in \inte_F(F_{n_0})$. By the homogeneity, there is a
compact Cantor manifold $K$ with respect to $D_{\mathcal K}^{k}$ for
some $k>n$ satisfying
\begin{enumerate}
\item $D_{\mathcal{K}}(K)=k+2$,
\item   $x\in K$,
\item $\diam K<\eta=d(x,F\setminus \inte_F (F_{n_0}))$.
\end{enumerate}
Since the set $K$ is not contained in $F$, we can assume that  there
is a point $a\in K\cap U$. Then, for
$0<\epsilon<\min\{\eta,d(a,X\setminus U)\}$  there is $\delta>0$ as
in Theorem~\ref{Effros}. Choosing a point $b\in V$, $d(x,b)<\delta$,
we obtain an $\epsilon$-homeomorphism $h:X\to X$ such that $h(x)=b$.
Then $h(K)$ is a Cantor manifold with respect to $D_{\mathcal
K}^{k}$ which is separated by a subset of $F_{n_0}$, a
contradiction.

\end{proof}

Next two propositions  are easy consequences of the definition of a
strong Cantor manifold.

\begin{pro}\label{p4}
Let $X$ be a  space satisfying condition~\eqref{eq1}
and let $K\subset X$ be a strong Cantor manifold with respect to an
admisible class $\mathcal C$. Then there exists exactly one $i$ such
that $K\subset F_i$.
\end{pro}

\begin{pro}\label{p5} Let $X$ be a locally compact Cantor manifold with respect to an
admissible class $\mathcal C$.   Assume  $X$ is not a strong Cantor
manifold with respect to $\mathcal C$ and no open non-empty subspace
of $X$ belongs to $\mathcal C$. Then $X$ satisfies
condition~\eqref{eq1}, i.e., $$X =\bigcup _{i=0}^{\infty} F_i \quad
\text{ with }\quad \bigcup _{i\ne j} (F_i \cap F_j)\in\mathcal C,$$
where  the sets  $F_i$ are proper, closed subsets of $X$ which
additionally satisfy
\begin{itemize}\item[(i)] no finite sum of $F_i$'s covers
$X$,\item[(ii)] $\inte F_i\ne\varnothing$ for each $i$,\item[(iii)]
$ F_i\cap \inte F_j=\varnothing$ for each $i\ne j$.
\end{itemize}
\end{pro}

\begin{proof}
Part (i) is a direct consequence of $X$ being a Cantor manifold with
respect to $\mathcal C$. To prove  (ii) and (iii) we can assume that
$\inte F_0\ne\varnothing$ by the Baire Category Theorem. Then, since
$F_0\ne X$, the open set $U_0=X\setminus \cl(\inte F_0)$ is
non-empty and  is contained in the union $(F_0\setminus \cl(\inte
F_0))\cup F_1\cup F_2 \cup \dots$, so there exits $n_1>n_0=0$ such
that $\inte F_{n_1}\ne\varnothing$. The open set
$U_1=X\setminus(\cl(\inte F_0)\cup\cl(\inte F_{n_1}))$ is non-empty
by (i) and it is contained in $(F_0\setminus \cl(\inte F_0))\cup
(F_{n_1}\setminus \cl(\inte F_{n_1}))\cup F_2\cup\dots$, etc. We
obtain a subsequence $n_0<n_1<n_2<\dots$ such that the sets
$F_{n_i}$ have non-empty interiors. Redefining $F'_0=F_0$,
$F'_i=F_{n_{i-1}+1}\cup\dots\cup F_{n_i}$, we get the representation
$X=\bigcup_{i=0}^\infty F'_i$ satisfying (ii). Notice that $\inte
F'_i\cap \inte F'_j=\varnothing$ if $i\ne j$. Indeed, otherwise this
intersection  would be a non-empty open subset of $X$, so it does
not belong to $\mathcal C$.  Since this open set is an
$F_\sigma$-subset of $F'_i\cap F'_j$, we get   $F'_i\cap F'_j\notin
\mathcal C$, a contradiction with $\bigcup_{i\ne j}(F'_i\cap
F'_j)\in\mathcal C$. Therefore, putting $F''_i=F'_i\setminus
\bigcup_{j\ne i}(\inte F'_j)$ we obtain the representation
$X=\bigcup_{i=1}^\infty F''_i$ with (i--iii) satisfied.
\end{proof}

\begin{thm}\label{thm2} Each  homogeneous continuum $X\notin\mathcal{C}$
is a strong Cantor manifold with respect to class $\mathcal{C}$
provided that:
 \begin{enumerate}
\item
$\mathcal{C}$ is any of the  following three classes:
$\mathcal{WID}$, $\mathbf C$, $\mathcal D_{\mathcal K}^{n-2}$
 $($in the latter case  we additionally assume
$D_{\mathcal{K}}(X)=n$$)$;\\
or
 \item
$\mathcal{C}=\mathcal D_{\mathcal K}^{<\infty}$  and   $X$ does not
contain closed subsets of arbitrary large finite dimension.
\end{enumerate}
\end{thm}

\begin{proof}
Suppose  $X$ is not a strong Cantor manifold with respect to
$\mathcal{C}$. Then, by Theorem~\ref{TC}, $X$ has a representation
$X=\bigcup_{i=0}^\infty F_i$ as in Proposition~\ref{p5}. By
Theorems~\ref{main1}, \ref{inf}, \ref{C}, \ref{wid} and
Proposition~\ref{mazur}
 $X$ contains  a  strong Cantor manifold with
respect to $\mathcal{C}$.
   By homogeneity, we can assume that any
point of $X$ belongs to such a strong Cantor manifold in $X$.

\begin{claim}\label{claim1}
If a strong Cantor manifold $K\subset X$ with respect to
$\mathcal{C}$ intersects $Y_i=\bd(\inte F_i)$, then $K\subset Y_i$.
 \end{claim}
Indeed, let $x\in K\cap\bd(\inte F_i)$ and suppose $K$ is not a
subset of $\bd(\inte F_i)$.  In the case where there exists $a\in
K\cap\inte F_i$, we can apply the Effros Theorem for
$0<\epsilon<d(a,(X\setminus\inte F_i)$ to find a $\delta>0$ such
that if a point $y\in X\setminus F_i$ is chosen with
$d(x,y)<\delta$, then there exists an $\epsilon$-homeomorphism
$h\colon X\to X$ that maps $x$ onto $y$. Then $h(K)$ is a strong
Cantor manifold with respect to $\mathcal{C}$ which intersects
$\inte F_i$ and another set $F_j$. This is however impossible by
Propositions~\ref{p4} and ~\ref{p5} (iii). In the case where there
is a point $a\in K\cap (X\setminus \cl(\inte F_i))$, we use an
$\epsilon$-homeomorphism $h$ which maps $x$ to a point in $\inte
F_i$ for $\epsilon<d(a,\cl(\inte F_i))$. The continuum  $h(K)$,
containing points in $\inte F_i$ and in $X\setminus \cl(\inte
F_i))$, must intersect $\bd(\inte F_i)$. Since $h(K)$ is a strong
Cantor manifold, we come to the former case above.

Let $K\subset Y_0=Y_X$ be a  strong Cantor manifold with respect to
$\mathcal{C}$. Define by transfinite induction:
\begin{align}\label{eq2} & K_0=K,\quad  K_{\alpha+1}=\cl\left(\bigcup\{h(K_\alpha):
h(K_\alpha)\cap K_\alpha\ne\varnothing, h\in H(X)\}\right)\\
&\text{and}\quad K_\alpha=\cl(\bigcup_{\beta<\alpha}
K_\beta)\quad\text{for  limit ordinals $\alpha$}\notag.\end{align}

There exists a countable ordinal $\gamma$ such that
$K_\gamma=K_{\gamma+1}=\dots$~\cite[Theorem 3, p. 258]{Ku}. Denote
$G_0=K_\gamma$.
\begin{claim}\label{claim2}
$G_0$  is a continuum contained in $Y_X$ and the group $H(X)$ acts
transitively on $G_0$.\end{claim} This follows  from~\eqref{eq2}, by
the homogeneity of $X$ and by Claim~\ref{claim1}.

\begin{claim}\label{claim3}$G_0$ is a strong Cantor manifold with respect to
$\mathcal{C}$.\end{claim} Suppose not. Then Claim~\ref{claim2}
allows us to repeat all the above considerations substituting $G_0$
for $X$ as the underlying space but  keeping the whole group $H(X)$
to act transitively on $G_0$. In particular, since $K\subset G_0$,
we get $K\subset Y_{G_0}\subsetneq G_0$ and  definition~\eqref{eq2}
gives $G_0\subset Y_{G_0}\subsetneq G_0$, a contradiction.
\begin{claim}\label{claim4}The collection $\mathcal G=\{h(G_0):h\in H(X)\}$
is a continuous decomposition of $X$.\end{claim} Observe that   each
two  distinct $G,G'\in\mathcal G$ are disjoint (see~\eqref{eq2}) and
if $h(G)\cap G'\ne\varnothing$, then $h(G)=G'$ for any $h\in H(X)$.
The continuity of the decomposition easily follows from the Effros
Theorem (cf.~\cite{Ro}).

To get a final contradiction, consider a correspondence $s:\mathcal
G\to \{F_1,F_2,\dots\}$ such that $G\subset s(G)$. By
Proposition~\ref{p4}, $s$ is a well defined function. Notice that
$$s^{-1}(F_i)\subset F_i,\quad\text{for each $i$,}\quad\text{and}\quad X=s^{-1}(F_1)\cup s^{-1}(F_2)\cup\dots$$
Since the decomposition $\mathcal G$ is continuous, the sets
$s^{-1}(F_i)$ are closed in $X$. It follows from the Sierpi\'nski
Theorem~\cite[p. 440]{re:89} that  the intersection $s^{-1}(F_i)\cap
s^{-1}(F_j)$ is nonempty for some $i\ne j$. Thus, the intersection
contains an element of $\mathcal G$ which is a strong Cantor
manifold with respect to $\mathcal{C}$, a contradiction with
Proposition~\ref{p4}.

\end{proof}

We do not know if one can omit, in general, the extra hypothesis in
Theorem~\ref{thm2}(2) that $X$ does not contain closed subsets of
arbitrary large finite dimension for $\mathcal{C}=\mathcal
D_{\mathcal K}^{<\infty}$.

The property $(\alpha)$ \label{alpha} of an $n$-dimensional space
$X$ (originally considered  by Hurewicz in~\cite{Hu} for the
covering dimension of separable spaces) means that any
$n$-dimensional closed subset of $X$ has the non-empty interior. It
is known that all topological $n$-manifolds have property $(\alpha)$
and it was observed in \cite{Se} that $n$-dimensional locally
compact, locally homogeneous ANR's also have this property for the
covering dimension.

It is proved in~\cite[Theorem C]{Se} that, for the covering
dimension,   an $n$-dimensional, locally compact, connected, locally
homogeneous ANR-space $X$ is a Cantor $n$-manifold. Actually, the
assumption in that proof that $X$ is an ANR reduces just to property
$(\alpha)$ and the reasoning is  applicable  to dimension
$D_{\mathcal K}$, so we have the following proposition.

\begin{pro}[\cite{Se}]\label{seidel}
If $X$ is a locally compact, connected, locally homogeneous  space
with property $($$\alpha$$)$ and $D_{\mathcal K}(X)=n$, then $X$ is
an $n$-dimensional  Cantor manifold with respect to $\mathcal
D_{\mathcal K}^{n-2}$.
\end{pro}

The following theorem generalizes this result.

\begin{thm}\label{thm3} Under the hypotheses of Proposition~$\ref{seidel}$, the space $X$
is a locally connected  strong Cantor manifold with respect to
$\mathcal D_{\mathcal K}^{n-2}$.
\end{thm}
\begin{proof}
By Remark~\ref{r}, there exist arbitrary small $n$-dimensional
compact subsets of $X$.
 Therefore $X$ contains
 arbitrary small compact, $n$-dimensional,  strong Cantor
manifolds with respect to $\mathcal D_{\mathcal K}^{n-2}$
(Theorem~~\ref{main1}). The local homogeneity and property
($\alpha$) guarantee that $X$ has a basis consisting of such strong
Cantor manifolds. In particular, $X$ is locally connected. Moreover,
by Proposition~\ref{seidel},  $X$ is an $n$-dimensional  Cantor
$n$-manifold with respect to $\mathcal D_{\mathcal K}^{n-2}$.
Suppose that $X$ is not a strong Cantor manifold with respect to
$\mathcal D_{\mathcal K}^{n-2}$. Then we can apply
Proposition~\ref{p5}. Since each point of $X$ is contained in the
interior of  an  $n$-dimensional  strong Cantor manifold $K$ and $K$
is contained in only one $F_i$ (see Proposition~\ref{p4}), it
follows that $F_i=\inte F_i$ for each $i$. This contradicts the
connectedness of $X$.
\end{proof}

\section{Remarks on property $(\alpha)$ for dimension $D_{\mathcal K}$}

We are going to propose an extension property $(H)$ which implies
property  $(\alpha)$ (see page \pageref{alpha}) for dimension
$D_{\mathcal K}$ in the class of  compact spaces. It is extracted
from a proof in~\cite{HW} and seems to be a natural and convenient
criterion  for  deriving property  $(\alpha)$ in many cases.

\begin{dfn} Let $\mathcal K$ be a given stratum. A space $X$ with an open
basis  $\mathcal U$ has property $(H)$ if
\begin{itemize}\item[$(H)$]
$D_{\mathcal K}(\bd U)\le n-1$  and any mapping $f:\bd U\to K_{n-1}$
extends over $(\cl U)\setminus V$ for each $U\in\mathcal U$ and any
open, nonempty subset  $V$ of $U$.
\end{itemize}
\end{dfn}
As in the case of property  $(\alpha)$, natural examples of spaces
with property $(H)$ are  manifolds (with or without boundaries).
Other examples include $n$-manifolds from which a  sequence (or
finite number of sequences) of mutually disjoint  open $n$-cells
converging to a point (or to different points, resp.) is removed. A
simple triod $T$ and the product $T\times \I$ have property
$(\alpha)$ but they do not have property $(H)$.

The proof of Theorem~\ref{thm1} follows the idea of~\cite[Theorem VI
12]{HW}. A key ingredient of the proof is the following lemma
(cf.~\cite[A), p.96]{HW}).

\begin{lem}\label{l1} Suppose a space $X$ has an open basis $\mathcal U$ satisfying $(H)$.
Let $Y$ be a closed subset of $X$. If $a\in Y$ is a boundary point
of $Y$ and $a\in U\in\mathcal U$, then any map of $Y\setminus U$
into $K_{n-1}$ can be extended over $Y$.
\end{lem}
 \begin{proof}
Let
$$\quad U'=U\cap Y\quad\text{and}\quad B=\bd U.$$
If $f\colon Y\setminus U'\to K_{n-1}$, then the restricted map
$f|_{(Y\setminus U')\cap B}$ extends to a map $f'$ over $B$, since
$(Y\setminus U')\cap B$ is a closed subset of $B$ and  $D_{\mathcal
K}(B)\le n-1$. Next, the map $f'$ can be extended to a map
$$g:(\cl
U)\setminus (U\setminus Y)\to K_{n-1}.$$
 Finally, the map
$\overline{f}: Y\to K_{n-1}$ given by
$$\overline{f}(x)=\begin{cases}g(x) &\text{for $x\in U',$}\\ f(x)
&\text{for $x\in Y\setminus U'$}\end{cases}$$ extends $f$.
\end{proof}

\begin{thm}\label{thm1}
If an $n$-dimensional compact space $X$ satisfies $(H)$, then $X$
has property $(\alpha)$.
\end{thm}
\begin{proof}
Let $Y$ be a closed $n$-dimensional subset of $X$. There exist a
closed subset $C$ of $Y$ and a map $f\colon C\to K_{n-1}$ which
cannot be extended over $Y$. By the compactness of $Y$ and Zorn
Lemma,  there exists a minimal (with respect to the inclusion)
closed subset $K$ of $Y$ such that $f$ is not extendable over $C\cup
K$. Then the set $K\setminus C$ is non-empty and we will  show that
it is open in $X$. Let $a\in K\setminus C$. Take $U\in\mathcal U$
such that $a\in U\subset \cl U\subset X\setminus C$. Since
$K\setminus U$ is a closed proper subset of  $K$, there exists an
extension $F\colon C\cup(K\setminus U)\to K_{n-1}$ of $f$. The  map
$F|_{\cl(K\setminus C)\setminus U}$ cannot be extended over
$\cl(K\setminus C)$ because if  $F'$ were such an extension,
then the map $G\colon C\cup K\to  K_{n-1}$ given by $$G(x)=\begin{cases}F'(x) &\text{if $x\in K\setminus C$,}\\
F(x) &\text{if   $x\in C$}\end{cases}$$ would extend $f$. By
Lemma~\ref{l1}, $a\in\inte(\cl(K\setminus C))$, hence
$a\in\inte(K\setminus C)\subset \inte Y$.
\end{proof}

\begin{qu} Are properties $(\alpha)$ and $(H)$ equivalent for finite-dimen\-sio\-nal
(locally) homogeneous compact spaces?
\end{qu}


\end{document}